\documentclass{amsart}

\usepackage{a4,xy,amssymb,diagrams}
\usepackage[hyperindex,pdfmark]{hyperref}

\xyoption{poly}
\xyoption{arc}
\xyoption{knot}
\xyoption{curve}
\xyoption{arrow}
\xyoption{all}

\makeindex

\newcommand{\dummy[1]}{{#1}}
\newcommand{\C}{\mathbb{C}}

\newcommand{\R}{\mathbb{R}}
\newcommand{\N}{\mathbb{N}}
\newcommand{\Z}{\mathbb{Z}}

\newcommand{\PP}{\mathbb{P}}

\newcommand{\Af}{\mathbb{A}}

\newcommand{\wis}[1]{{\text{\em \usefont{OT1}{pag}{m}{n} #1}}}

\newtheorem{theorem}{Theorem}[section]
\newtheorem{lemma}[theorem]{Lemma}
\newtheorem{proposition}[theorem]{Proposition}

\theoremstyle{definition}

\theoremstyle{remark}

\begin{document}

%%%%%%%%%%%%%%%%%%%%%
%%%some commands
%%%%%%%%%%%%%%%%%%%%%
\sloppy
%%%%%%%%%%%%%%%%%%%%%
%%%end of commands
%%%%%%%%%%%%%%%%%%%%%

\title{Local quivers and stable representations}

% Remove or comment out any unused author tags.
\author{Jan Adriaenssens}
\address{Universiteit Antwerpen (RUCA) \\ B-2020 Antwerp (Belgium)}
\email{adri@ruca.ua.ac.be}
\urladdr{http://math.ruca.ua.ac.be/$\sim$adri/}
\author{Lieven Le Bruyn}% first author
\address{Universiteit Antwerpen (UIA) \\ B-2610 Antwerp (Belgium)}
%\curraddr{}
\email{lebruyn@uia.ua.ac.be}
\urladdr{http://win-www.uia.ac.be/u/lebruyn/}
\thanks{Research Director of the FWO (Belgium)}

\begin{abstract}
In this paper we introduce and study the local quiver as a tool to investigate the
\'etale local structure of moduli spaces of $\theta$-semistable representations of
quivers. As an application we determine the dimension vectors associated to irreducible
representations of the torus knot groups $G_{p,q} = \langle~a,b~\mid~a^p=b^q~\rangle$.
\end{abstract}

\maketitle

\section{Introduction.}

Polynomial invariants of quiver representations are generated by taking traces along oriented
cycles in the quiver, see \cite{LeBruynProcesi:1990}. Hence, if we have a quiver $Q$ without oriented
cycles on $k$ vertices, the quotient variety of $\alpha$-dimensional representations $rep_{\alpha}~Q$ under the
basechange group $GL(\alpha)$ consists of one point corresponding to the 
unique $\alpha$-dimensional semi-simple representation which is the direct sum of the vertex
spaces. Still, there may be lots of semi-invariants on $rep_{\alpha}~Q$ with respect to a
character $\chi = \prod~det(g_i)^{t_i}$ of $GL(\alpha)$, determined by an integral vector
$\theta = (t_1,\hdots,t_k) \in \Z^k$. If we collect all semi-invariants with weight $n \chi$ for
all $n \in \N$ we obtain a positively graded connected algebra with associated projective
variety $M^{ss}_{\alpha}(Q,\theta)$ the moduli space of $\theta$-semistable $\alpha$-dimensional
representations of $Q$, see \cite{King:1994}. The points of $M^{ss}_{\alpha}(Q,\theta)$ correspond
to isomorphism classes of direct sums of $\theta$-stable representations of $Q$. Recall that an
$\alpha$-dimensional representation $V$ is $\theta$-stable (resp. $\theta$-semistable) if for all subrepresentation $W \rInto V$
we have $\theta(W) = \sum_i t_ib_i > 0$ (resp. $\geq 0$) where $\beta = (b_1,\hdots,b_k)$ is the dimension vector
of $W$. Not much is known on the geometry of these moduli spaces in particular of their
singularities.

As an illustrative example, consider the extended $n$-arrow Kronecker quiver $K_n$
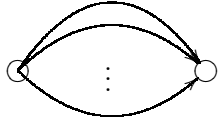
\[
\xy 0
\POS (25,0) *\cir<4pt>{}*+{} ="v2"
   , (0,0) *\cir<4pt>{}*+{} ="v3", (12,0) *+{\vdots}
\POS"v3" \ar@/^4ex/ "v2"
\POS"v3" \ar@/^6ex/ "v2"
\POS"v3" \ar@/_4ex/ "v2"
\endxy
\]
Take the dimension vector $\alpha = (1,r)$ with $r \leq n$ and the character on $GL(\alpha) = \C^* \times GL_r$
determined by $\theta = (-r,1)$. An $\alpha$-dimensional representation $V \in rep_{\alpha}~K_n$ is
fully determined by the $r \times n$ matrix whose $i$-th column is the linear map corresponding to
the $i$-th arrow. Clearly, a $\theta$-stable representation $V \in rep_{\alpha}~K_n$ must have an
$r \times n$ matrix of maximal rank for if the rank were $s < r$, then $V$ has a subrepresentation
of dimension vector $\beta = (1,s)$ and $\theta(W) = -r.1+1.s < 0$. The basechange action is given
by left multiplication by $GL_r$ and scalar multiplication by $\lambda^{-1}$ for $\lambda \in \C^*$.
Therefore, in this case we have that the moduli space of $\theta$-stable quiver representations
\[
M^{ss}_{\alpha}(K_n,\theta) \simeq Grass(r,n) \]
the Grassmannian of $r$-dimensional subspaces in $\C^n$. Smoothness of this moduli space follows
from the fact that $\alpha$ is minimal among the dimension vectors having $\theta$-stable
representations. Things become more complicated for multiples $k \alpha$ 
because in this case the moduli spaced $M^{ss}_{k \alpha}(K_n,\theta)$ contains points corresponding to
direct sums of $\theta$-stable representation such as
\[
V = W^{\oplus k} \quad \text{for $W \in rep_{\alpha}~K_n$ $\theta$-stable.} \]
We will show in this paper that the \'etale local structure of $M^{ss}_{k \alpha}(K_n,\theta)$
near the point determined by $V$ is that of the quotient variety of $r(n-r)$-tuples of $k \times
k$ matrices under simultaneous conjugation near the point determined by the tuple of zero matrices.
It is well known that this quotient variety is smooth only if $k = 1$ or $r(n-r) = 2 = k$. Therefore,
the only moduli spaces $M^{ss}_{k \alpha}(K_n,\theta)$ which are smooth are the Grassmannians
described before and
\[
M^{ss}_{(2,2)}(K_3,\theta)~\text{for $\theta = (-1,1)$}~\quad \text{and} \quad 
M^{ss}_{(2,4)}(K_3,\theta)~\text{for $\theta = (-2,2)$} \]
Both these varieties are isomorphic to $M_{\PP^2}(0,2) \simeq \PP^5$, the moduli space of semistable rank $2$
bundles on $\PP^2$ with Chern classes $c_1 = 0$ and $c_2 = 2$ by results of Barth and Hulek
\cite{Barth:1977},
\cite{Hulek:1979} (and reflexion functors). 

The general result is as follows. Let
$\xi \in M^{ss}_{\alpha}(Q,\theta)$ be a point corresponding to a representation
\[
V_{\xi} = W_1^{\oplus m_1} \oplus \hdots \oplus W_l^{\oplus m_l} \]
where $W_i$ is a $\theta$-stable representation of $Q$ of dimension vector $\beta_i$. Consider the
local quiver $Q'$ on $l$ vertices $w_1,\hdots,w_l$ (corresponding to the distinct $\theta$-stable
components) with
$
\delta_{ij}-\chi_Q(\beta_i,\beta_j)$ arrows from $w_i$ to $w_j$ where $\chi_Q$ is the Euler form of the
quiver $Q$. Consider the dimension vector $\alpha' = (m_1,\hdots,m_l)$ for $Q'$ determined by the
multiplicities of the $\theta$-stable components.

\begin{theorem} With notations as before, there is an \'etale isomorphism between a neighborhood of
$\xi$ in the moduli space $M^{ss}_{\alpha}(Q,\theta)$ and a neighborhood of the trivial representation
$\overline{0}$ in the quotient variety $iss_{\alpha'}~Q'$ of semi-simple $\alpha'$-dimensional
representations of the local quiver $Q'$.
\end{theorem}

The proof uses two ingredients : first we use the theory of universal localizations to reduce the
Zariski local structure of $M^{ss}_{\alpha}(Q,\theta)$ near $\xi$ to that of $\alpha$-dimensional
semi-simple representations of a universal localization $\C Q_{\sigma}$ of the path algebra. Next,
we apply the Luna slice machinery to this quotient variety to obtain the result.

The theorem can also be used to obtain an inductive method to determine all dimension vectors
of $\theta$-stable representations.

\begin{theorem} With notations as above, $\alpha$ is the dimension vector of a $\theta$-stable
representation if and only if $\alpha'$ is the dimension vector of a simple representation of the
local quiver $Q'$.
\end{theorem}

The utility of this result follows from the easy description of the dimension vectors of simple
quiver representations given in \cite{LeBruynProcesi:1990} by a set of inequalities, rather than an
inductive procedure. Finally, we apply this result to determine the dimension vectors associated to
irreducible representations of the torus knot groups, among which is the braid group on three
strings.

\section{Moduli spaces of quiver representations.}

In this section we recall briefly the basic setting.
Recall that a {\em quiver} $Q$ is a finite directed graph on a set of
vertices
$Q_v = \{ v_1,\hdots,v_k \}$, having a finite set $Q_a = \{ a_1,\hdots,a_l
\}$ of arrows,
where we allow loops as well as multiple arrows between vertices. An arrow $a$
with starting vertex
$s(a) = v_i$ and terminating vertex $t(a) = v_j$ will be depicted as
$\xy
\POS (15,0) *\cir<4pt>{}*+{\txt\tiny{i}} ="v2"
   , (0,0) *\cir<4pt>{}*+{\txt\tiny{j}} ="v3"
\POS"v2" \ar "v3"_{a} \endxy$. The quiver information is encoded in the
{\em Euler form}
which is the bilinear form on $\Z^k$ determined by the matrix $\chi_Q \in
M_k(\Z)$ with
\[
\chi_{ij} = \delta_{ij} - \#~\{~a \in Q_a~\mid~\xy
\POS (15,0) *\cir<4pt>{}*+{\txt\tiny{i}} ="v2"
   , (0,0) *\cir<4pt>{}*+{\txt\tiny{j}} ="v3"
\POS"v2" \ar "v3"_{a} \endxy~\}
\]
The {\em path algebra} $\C Q$ of a quiver $Q$ has as basis the set of all
oriented paths
$p = a_{i_u} \hdots a_{i_1}$ of length $u \geq 1$ in the quiver, that is
$s(a_{i_{j+1}}) = t(a_{i_j})$
together with the vertex-idempotents $e_i$ of length zero. Multiplication in
$\C Q$ is induced by (left) concatenation of paths. More precisely, $1= e_1
+ \hdots + e_k$ is a
decomposition of $1$ into mutually orthogonal idempotents and further we define
\begin{itemize}
\item{$e_j.a$ is always zero unless \ \  $\xy
\POS (0,0) *\cir<4pt>{}*+{\txt\tiny{j}} ="v2"
   , (15,0) *\cir<4pt>{} ="v3"
\POS"v3" \ar "v2"_{a} \endxy$ \ \ in which case it is the path $a$,}
\item{$a.e_i$ is always zero unless \ \ $\xy
\POS (15,0) *\cir<4pt>{}*+{\txt\tiny{i}} ="v2"
   , (0,0) *\cir<4pt>{} ="v3"
\POS"v2" \ar "v3"_{a} \endxy$ \ \ in which case it is the path $a$,}
\item{$a_i.a_j$ is always zero unless \ \ $\xy
\POS (15,0) *\cir<4pt>{} ="v2"
   , (0,0) *\cir<4pt>{} ="v1"
   , (30,0) *\cir<4pt>{} ="v3"
\POS"v2" \ar "v1"_{a_i}
\POS"v3" \ar "v2"_{a_j} \endxy$ \ \ in which case it is the path $a_ia_j$.}
\end{itemize}
Path algebras of quivers are the archetypical examples of {\em formally
smooth algebras} as
introduced and studied in \cite{CuntzQuillen:1995}.

For a dimension vector $\alpha = (a_1,\hdots,a_k) \in \N^k$ we define the
representation space
\[
rep_{\alpha}~Q = \underset{\xy
\POS (15,0) *\cir<4pt>{}*+{\txt\tiny{i}} ="v2"
   , (0,0) *\cir<4pt>{}*+{\txt\tiny{j}} ="v3"
\POS"v2" \ar "v3" \endxy}{\bigoplus}~M_{a_j \times a_i}(\C) \]
which is an affine space with natural action of the basechange group
$GL(\alpha) = GL_{a_1} \times \hdots \times GL_{a_k}$. The points of $rep_{\alpha}~Q$
correspond to left $\C Q$-modules $V$ such that $dim_{\C}~e_i.V = a_i$. The $GL(\alpha)$-orbits
correspond to isomorphism classes as $\C Q$-modules. By abuse of notation we speak of
morphisms (in particular subrepresentations) and extensions of representations for the
corresponding module theoretic notions.

Let $\theta = (t_1,\hdots,t_k) \in \Z^k$ and define for a $Q$-representation $W$ of dimension vector
$\beta = (b_1,\hdots,b_k)$ the number $\theta(W) = t_1b_1 + \hdots + t_k b_k$.
Following A. King in \cite{King:1994} we say that a representation $V \in rep_{\alpha}~Q$ is
$\theta$-semistable if $\theta(V) = 0$ and for every subrepresentation $W \subset V$ we have
$\theta(W) \geq 0$. A $\theta$-semistable representation $V$ is said to be $\theta$-stable if
the only subrepresentations $W \subset V$ with $\theta(W) = 0$ are $0$ and $V$ itself.

$\theta$ defines a character on $GL(\alpha)$ by sending a $k$-tuple $(g_1,\hdots,g_k) \in GL(\alpha)$
to $det(g_1)^{t_1} \hdots det(g_k)^{t_k} \in \C^*$. With $SI(Q,\theta)$ we denote the space of all
semi-invariant functions $f$ on $rep_{\alpha}~Q$ of weight $\theta$, that is $f \in \C[rep_{\alpha}~Q]$ such that
\[
f(g.V) = \theta(g) f(V) \quad \text{for all $g \in GL(\alpha)$ and $V \in rep_{\alpha}~Q$.}
\]
In \cite[Prop.3.1]{King:1994} it is shown that $V \in rep_{\alpha}~Q$ is $\theta$-semistable if
and only if there is a semi-invariant $f \in SI(Q,n \theta)$ for some $n \geq 0$ such that
$f(V) \not= 0$. For this reason we consider the positively graded algebra
\[
R_{\theta} = \oplus_{n=0}^{\infty} SI(Q,n \theta) \]
with part of degree zero $SI(Q,0) = \C[rep_{\alpha}~Q]^{GL(\alpha)}$ the ring of polynomial
quiver-invariants which is, by a result of \cite{LeBruynProcesi:1990} generated by traces along
oriented cycles in the quiver $Q$. In particular,if $Q$ is a quiver without oriented cycles,
then $R_{\theta}(0) = \C$ and we obtain a projective variety
\[
\wis{proj}~R_{\theta} = M^{ss}_{\theta}(Q,\alpha) \]
called the moduli space of $\theta$-semistable representations of the quiver $Q$ of dimension
vector $\alpha$. 

By \cite[Prop. 3.2]{King:1994} a point $\xi \in M_{\theta}^{ss}(Q,\alpha)$ determines uniquely the
orbit of a
representation $V_{\xi} \in rep_{\alpha}~Q$ such that
\[
V = W_1^{\oplus m_1} \oplus \hdots W_l^{\oplus m_l} \]
where the $W_i$ are $\theta$-stable representations of $Q$ of dimension vector $\beta_i$ which
occur in $V_{\xi}$ with multiplicity $m_i$. In particular, the {\em semistable representation type}
$\tau = (m_1,\beta_1;\hdots;m_l,\beta_l)$ of $\xi$ or $V_{\xi}$ satisfies
\[
\alpha = m_1 \beta_1 + \hdots + m_l \beta_l \]
In this paper we will introduce a local quiver $Q_{\xi}$ (depending only on the semistable
representation type) as a tool to investigate the \'etale (that is, analytic) local structure
of the moduli space $M_{\theta}^{ss}(Q,\theta)$ in a neighborhood of $\xi$.

\section{Universal localizations.}

In this section we will show that one can cover the projective moduli space $M_{\theta}^{ss}(Q,\alpha)$
by affine open sets, each of which isomorphic to the quotient variety 
of $rep_{\alpha}~\C Q_{\sigma}$ where $\C Q_{\sigma}$
is a certain universal localization of the path algebra $\C Q$. First, we need to recall a result
on the generation of semi-invariants of quiver representations, see for example \cite{SchofieldVdB:1999},
\cite{DerksenWeyman:1999} and \cite{DomokosZubkov}. Throughout, we will assume that $Q$ is a quiver
without oriented cycles.

Fix a character $\theta = (t_1,\hdots,t_k)$ and divide the set of vertex-indices into a left
set $L = \{ i_1,\hdots,i_u \}$ consisting of those $1 \leq i \leq k$ such that $t_i \leq 0$ and
a right set $R = \{ j_1,\hdots,j_v \}$ consisting of those $1 \leq j \leq k$ such that
$t_j \geq 0$ (observe that $L \cap R$ may be non-empty).

For every vertex $v_i$ we define the indecomposable projective module $P_i = \C Q e_i$ spanned by
all paths in the quiver $Q$ starting at $v_i$. As a consequence we have that
$Hom_Q(P_i,P_j)$ is spanned by all paths $[j,i]$ in the quiver $Q$ from vertex $v_j$ to vertex $i$.
Because $Q$ has no oriented cycles, these spaces are finite dimensional.

For a fixed integer $n$ we consider the set $\Sigma_{\theta}(n)$ of all $\C Q$-module morphisms
\[
\begin{diagram}
P_{i_1}^{\oplus -n t_{i_1}} \oplus \hdots \oplus P_{i_u}^{\oplus - n t_{i_u}} & \rTo^{\sigma} &
P_{j_1}^{\oplus n t_{j_1}} \oplus \hdots \oplus P_{j_v}^{\oplus  n t_{j_v}} \\
\Arrowvert_{\text{notation}} & & \Arrowvert_{\text{notation}} \\
P_{c_1} \oplus \hdots \oplus P_{c_p} & \rTo^{\sigma} & P_{d_1} \oplus \hdots \oplus P_{d_q}
\end{diagram}
\]
By the remark above, $\sigma$ can be described by an $(p = n \sum t_{i_l}) \times (q= n \sum t_{j_m})$ matrix 
$M_{\sigma}$ all
entries of which are linear combinations $p_{lm}$ of paths in the quiver $Q$ from vertex
$v_{d_m}$ to vertex $v_{c_l}$.
For $V \in rep_{\alpha}~Q$ we can substitute the arrow matrices $V(a)$ in the definition of
$p_{lm}$ and obtain a square matrix of size
$a_{c_l} \times a_{d_m}$. If we do this for every entry of $\sigma$ we obtain a
square matrix as $\theta(V) = 0$ which we denote by $\sigma(V)$. But then, the function
\[
d_{\sigma}(V) = det~\sigma(V)~:~rep_{\alpha}~Q \rTo \C \]
is a semi-invariant of weight $n \theta$. The semi-invariants $SI(Q,n \theta)$ are spanned by
such determinantal semi-invariants.

The universal localization $\C Q_{\sigma}$ of the path algebra $\C Q$ with respect to the
morphism $\sigma$ is the algebra extension $\C Q \rTo^{j_{\sigma}} \C Q_{\sigma}$ such that
the extended morphism $\C Q_{\sigma} \otimes_{\C Q} \sigma$ becomes invertible and $j_{\sigma}$ has
the universal property with respect to this condition.

It is easy to see that $\C Q_{\sigma}$ is an affine $\C$-algebra. For, add an additional
$q \times p$ arrows $n_{lm}$ to the quiver $Q$ where $n_{lm}$ is an arrow from vertex $v_{c_l}$
to vertex $v_{d_m}$ and consider the $q \times p$ matrix 
\[
N_{\sigma} = \begin{bmatrix} n_{11} & \hdots & n_{1p} \\
\vdots & & \vdots \\
n_{q1} & \hdots & n_{qp} \end{bmatrix} \]
Then, the universal localization $\C Q_{\sigma}$ is the quotient of the path algebra of this
extended quiver modulo the ideal of relations given by the entries of the following matrix
identities
\[
M_{\sigma}.N_{\sigma} = \begin{bmatrix} v_{c_1} & & 0 \\
& \ddots & \\
0 & & v_{c_p} \end{bmatrix} \qquad
N_{\sigma}.M_{\sigma} = \begin{bmatrix}
v_{d_1} & & 0 \\
& \ddots & \\
0 & & v_{d_q} \end{bmatrix} \]
It follows from this description that for a dimension vector $\beta = (b_1,\hdots,d_k)$ the
variety of $\beta$-dimensional representations of $\C Q_{\sigma}$ is the open (but possibly
empty) subvariety of $rep_{\beta}~Q$ consisting of those representations $W$ such that
$\sigma(W)$ is a square invertible matrix. In particular, it follows that 
$rep_{\beta}~\C Q_{\sigma} = \emptyset$ unless $\theta(\beta) = \sum t_i b_i = 0$.

Returning to the study of $\theta$-semistable quiver representations and their moduli
spaces, we define
\[
(~rep_{\alpha}~\C Q_{\sigma} = )~X_{\sigma}(\alpha) = \{ V \in rep_{\alpha}~Q~\mid~d_{\sigma}(V) \not= 0 \}
\]
Because $d_{\sigma}$ is a semi-invariant of weight $n \theta$ it follows that $X_{\sigma}(\alpha)$
consists of $\theta$-semistable representations.

\begin{lemma} For $V \in X_{\sigma}(\alpha)$, the following are equivalent 
\begin{enumerate}
\item{$V$ is a $\theta$-stable representation.}
\item{$\C Q_{\sigma} \otimes V$ is a simple $\alpha$-dimensional representation of $\C Q_{\sigma}$.}
\end{enumerate}
\end{lemma}

\begin{proof}
Let $V$ be $\theta$-stable and assume that $W \subset V$ is a proper sub $\C Q_{\sigma}$-module
of dimension vector $\beta$. Restricting $W$ to a representation of $Q$ we see that $W \in
rep_{\beta}~Q$ and is a subrepresentation of $V$. Because $W \in rep_{\beta}~\C Q_{\sigma}$ we
have $\theta(\beta) = \theta(W) = 0$, which is impossible as $V$ is $\theta$-stable.

Let $\C Q_{\sigma} \otimes V$ be a simple $\C Q_{\sigma}$ representation and assume that $W \subset V$
is a proper subrepresentation of dimension vector $\beta = (b_1,\hdots,b_k)$
with $\theta(W) \leq 0$. If $\theta(W) < 0$ then $-\sum n t_{i_l} b_{i_l} > \sum n t_{j_m} b_{j_m}$
whence $\sigma(W)$ has a kernel but this contradicts the fact that $\sigma(V)$ is invertible.
Hence, $\theta(W) = 0$ but then $\C Q_{\sigma} \otimes W$ is a proper subrepresentation of 
$\C Q_{\sigma} \otimes V$ contradicting simplicity.
\end{proof}

Let $rep^{\theta-ss}_{\alpha}~Q$ be the open set of $\theta$-semistable
representations in $rep_{\alpha}~Q$. Let $\xi \in M^{ss}_{\theta}(Q,\alpha)$ of semistable representation type
$\tau = (m_1,\beta_1;\hdots;m_l,\beta_l)$, then the fiber of the quotient map
\[
rep^{\theta-ss}_{\alpha}~Q \rOnto^{\pi_{\theta}} M^{ss}_{\theta}(Q,\alpha) \]
contains a unique closed orbit corresponding to
\[
V_{\xi} = W_1^{\oplus m_1} \oplus \hdots \oplus W_l^{\oplus m_l} \]
with the $W_i$ all $\theta$-stable. There is a determinantal semi-invariant $d_{\sigma} \in SI(Q,n \theta)$
for some $n$ such that $d_{\sigma}(V_{\xi}) \not= 0$, hence $V_{\xi}$ belongs to the affine open subset
$X_{\sigma}(\alpha)$ of $rep_{\alpha}~Q$ (note that $X_{\sigma}(\alpha)$ is also an affine open
$GL(\alpha)$-invariant subset of $rep^{\theta-ss}_{\alpha}~Q$). Every $\theta$-semistable representation in the fiber
$\pi_{\theta}^{-1}(\xi)$ has a filtration with factors the $W_i$. After tensoring with the 
universal localization $\C Q_{\sigma}$ we have that the $W_i$ become simple representations, 
$V_{\xi}$ is a semi-simple representation and the above mentioned filtration becomes a Jordan-H\"older
filtration. 

\begin{theorem}
The quotient map $\pi_{\theta}$ is locally isomorphic to the
quotient map
\[
\begin{diagram}
rep_{\alpha}~\C Q_{\sigma} = X_{\sigma}(\alpha) & \rInto & rep^{\theta-ss}_{\alpha}~Q \\
\dOnto^{\pi_{\sigma}} & & \dOnto^{\pi_{\theta}} \\
iss_{\alpha}~\C Q_{\sigma} = X_{\sigma}(\alpha)/GL(\alpha) & \rInto & M_{\theta}^{ss}(Q,\alpha)
\end{diagram}
\]
assigning to an $\alpha$-dimensional $\C Q_{\sigma}$-module its semi-simplification, that is, the
direct sum of its Jordan-H\"older components. Because the affine open sets $X_{\sigma}(\alpha)$
cover $rep_{\alpha}^{\theta-ss}~Q$, the moduli space of $\theta$-semistable quiver representations
$M^{ss}_{\theta}(Q,\alpha)$ is locally isomorphic to quotient varieties $iss_{\alpha}~\C Q_{\sigma}$
for specific universal localizations $\C Q_{\sigma}$ of the path algebra.
\end{theorem}

\section{The local quiver.}

In this section we will prove that the local structure of the moduli space $M_{\theta}^{ss}(Q,\alpha)$
is determined by a new quiver setting which we will now describe. Let $\xi \in M_{\theta}^{ss}(Q,\alpha)$
of semistable representation type $\tau = (m_1,\beta_1;\hdots;m_l,\beta_l)$, that is the unique closed
orbit lying in the fiber $\pi^{-1}_{\theta}(\xi)$ is the isomorphism class of a direct sum
\[
V_{\xi} = W_1^{\oplus m_1} \oplus \hdots \oplus W_l^{\oplus m_l} \]
with $W_i$ a $\theta$-stable representation of dimension vector $\beta_i$. We introduce a new
quiver setting $(Q_{\xi},\alpha_{\xi})$ as follows :
\begin{itemize}
\item{$Q_{\xi}$ has $l$ vertices $w_1,\hdots,w_l$ corresponding to the distinct $\theta$-stable
components of $V_{\xi}$, and}
\item{the number of arrows from $w_i$ to $w_j$ is equal to
\[
\delta_{ij} - \chi_Q(\beta_i,\beta_j) \]
where $\chi_Q$ is the Euler form of the original quiver $Q$.}
\item{the dimension vector $\alpha_{\xi} = (m_1,\hdots,m_l)$ gives the multiplicities of
the stable summands.}
\end{itemize}
Observe that the local quiver setting depends only on the semistable representation type.

\begin{theorem} There is an \'etale isomorphism between an affine neighborhood of $\xi$ in the
moduli space $M_{\theta}^{ss}(Q,\alpha)$ and an affine neighborhood of the image $\overline{0}$ of the
zero representation in the quotient variety $iss_{\alpha_{\xi}}~Q_{\xi} = rep_{\alpha_{\xi}}~Q_{\xi}/GL(\alpha_{\xi})$
corresponding to the local quiver setting $(Q_{\xi},\alpha_{\xi})$.
\end{theorem}

Let $\xi \in M_{\theta}^{ss}(Q,\alpha)$ with corresponding $V_{\xi}$ having a decomposition into
$\theta$-stable representations $W_i$ as above. By the results of the preceding section we may
assume that $V_{\xi} \in X_{\sigma}(\alpha)$ where $X_{\sigma}(\alpha)$ is the affine
$GL(\alpha)$-invariant open subvariety of $rep_{\alpha}^{ss}~Q$ defined by the determinantal
semi-invariant $d_{\sigma}$.

We have seen that $X_{\sigma}(\alpha) \simeq rep_{\alpha}~\C Q_{\sigma}$ the variety of $\alpha$-dimensional
representations of the universal localization $\C Q_{\sigma}$. Moreover, if we define
$V_{\xi}' = \C Q_{\sigma} \otimes V_{\xi}$ and $W_i' = \C Q_{\sigma} \otimes W_i$ we have seen that
\[
V_{\xi}' = W_1^{'\oplus m_1} \oplus \hdots \oplus W_l^{'\oplus m_l} \]
is a decomposition of the semisimple $\C Q_{\sigma}$ representation $V'_{\xi}$ into its simple
components $W'_i$. Restricting to the affine smooth variety $X_{\sigma}(\alpha)$ we are in a
situation to apply the strong form of the Luna slice theorem, \cite{Luna:1973} or \cite{Slodowy:slice},
which we now recall.

Let $X$ be a smooth affine variety on which a reductive group $G$ acts. Let $x \in X$ be a point in
a closed $G$-orbit, then the stabilizer subgroup $G_x = \{ g \in G~\mid~g.x = x \}$ is again a
reductive group. Let $N_x$ be the normal space to the orbit $G.x$ of $x$, that is,
$N_x = \tfrac{T_x~X}{T_x~G.x}$, then $N_x$ is a finite dimensional $G_x$-representation. The Luna
slice theorem asserts the existence of a locally closed subvariety (the slice) $S \rInto X$
satisfying the following properties :
\begin{itemize}
\item{$S$ is a smooth affine variety with a $G_x$-action.}
\item{The action map $G \times S \rTo X$ sending $(g,s)$ to $g.s$ induces a $G$-equivariant
\'etale morphism from the associated fiber bundle
\[
G \times^{G_x} S \rTo^{\psi} X \]
with an affine image and such that the induced quotient map
\[
(G \times^{G_x} S)/G_x \simeq N_x/G_x \rTo^{\psi/G} X/G \]
is \'etale.}
\item{There is a $G_x$-equivariant morphism $S \rTo^{\phi} T_x~S = N_x$ with $\phi(x) = 0$,having an affine image
and such that the induced quotient map
\[
S/G_x \rTo^{\phi/G_x} N_x/G_x \]
is an \'etale map.}
\end{itemize}
That is, we have a commutative diagram of varieties
\[
\begin{diagram}
G \times^{G_x} N_x & & \lTo^{G \times^{G_x} \phi} & & G \times^{G_x} S & & \rTo^{\psi} & & X \\
\dOnto & & & & \dOnto & & & & \dOnto^{\pi} \\
N_x/G_x & & \lTo^{\phi/G_x} & & S / G_x & & \rTo^{\psi/G} & & X/G
\end{diagram}
\]
where the upper horizontal maps are all $G$-equivariant \'etale maps and the lower horizontal
maps are \'etale. If $\overline{x} = \pi(x)$, we deduce that the \'etale local structure of the
quotient variety $X/G$ in a neighborhood of $\overline{x}$ is \'etale isomorphic to the
local structure of the quotient variety $N_x/G_x$ in a neighborhood of $\overline{0}$ the point
corresponding to the fixed $G_x$-orbit of the zero vector in $N_x$.

We will now make this general result explicit in the case of interest to us. $rep_{\alpha}~\C Q_{\sigma}$
is a smooth affine variety with action of $GL(\alpha)$ by basechange. The point $V'_{\xi}$ lies in a
closed $GL(\alpha)$-orbit as closed orbits correspond to semisimple representations by the
result of Artin \cite{Artin:1969} and Voigt \cite{Gabriel:1974}. By a result of Gabriel's
\cite{Gabriel:1974} we know that the normal space to the orbit $GL(\alpha).V'_{\xi}$ can be
identified with the self-extensions
\[
N_{V'_{\xi}} = \dfrac{T_{V'_{\xi}}~rep_{\alpha}~\C Q_{\sigma}}{T_{V'_{\xi}}~GL(\alpha).V_{\xi}} =
Ext^1_{\C Q_{\sigma}}(V'_{\xi},V'_{\xi}) = \oplus_{i,j=1}^l~Ext^1_{\C Q_{\sigma}}(W'_i,W'_j)^{\oplus m_im_j} \]
By Schur's lemma we know that the stabilizer subgroup of the semisimple module $V'_{\xi}$ is equal
to
\[
G_{V'_{\xi}} \simeq GL(\alpha_{\xi}) = GL_{m_1} \times \hdots \times GL_{m_l} \rInto GL(\alpha) \]
and if we write out the action of this group on the self extensions we observe that it coincides
with the action of the basechangegroup $GL(\alpha_{\xi})$ on the representation space
$rep_{\alpha_{\xi}}~\Gamma$ of a quiver $\Gamma$ on $l$ vertices such that the number of arrows
from vertex $w_i$ to vertex $w_j$ is equal to the dimension of the extension group
\[
Ext^1_{\C Q_{\sigma}}(W'_i,W'_j) \]
Remains to prove that the quiver $\Gamma$ is our local quiver $Q_{\xi}$. For this we apply a
general homological result valid for universal localizations, \cite[Thm 4.7]{Schofield:1986}. If
$A_{\sigma}$ is a universal localization of an algebra $A$, then the category of left $A_{\sigma}$-modules
is closed under extensions in the category of left $A$-modules. Therefore,
\[
Ext^1_{A_{\sigma}}(M,N) = Ext^1_A(M,N) \]
for $A_{\sigma}$-modules $M$ and $N$. We apply this general result to get
\[
Ext^1_{\C Q_{\sigma}}(W'_i,W'_j) = Ext^1_Q(W_i,W_j) \]
Further, as the $W_i$ are $\theta$-stable representations of the quiver $Q$ we know that
$Hom_Q(W_i,W_j) = \delta_{ij} \C$. Finally, we use the homological interpretation of the Euler
form valid for all representations of quivers, see for example \cite{KraftRiedtmann}
\[
\chi_Q(\beta_i,\beta_j) = dim_{\C}~Hom_Q(W_i,W_j) - dim_{\C}~Ext^1_Q(W_i,W_j) \]
to deduce that $\Gamma = Q_{\xi}$, which finishes the proof of the theorem.

\par \vskip 4mm

An important application of this local description is to determine the smooth locus of the
moduli spaces $M_{\theta}^{ss}(Q,\alpha)$.

\begin{proposition} $\xi \in M_{\theta}^{ss}(Q,\alpha)$ is a smooth point if and only if the
ring of polynomial quiver invariants of the corresponding local quiver setting
\[
\C [rep_{\alpha_{\xi}}~Q_{\xi}]^{GL(\alpha_{\xi})} \]
is a polynomial ring.
\end{proposition}

As we know from \cite{LeBruynProcesi:1990} that the ring of invariants is generated by taking
traces along oriented cycles in the quiver $Q_{\xi}$ this criterium is rather effective. The
classification of all coregular quiver settings $(Q,\alpha)$ (that is such that the ring of polynomial
invariants is a polynomial algebra) seems to be a difficult problem though recently 
significant progress has been made by R. Bocklandt \cite{Bocklandt:2000}.

\section{$\theta$-stable representations.}

In this section we will give an algorithm to determine the set of dimension vectors
of $\theta$-stable representations of a quiver $Q$. In the following sections we will give an application
of this general result.

An elegant, though inefficient, algorithm to determine the dimension vectors of
$\theta$-(semi)stable representations is based on the general quiver representation results
of A. Schofield \cite{Schofield:1992}. If $\alpha$ is the dimension
vector of a $\theta$-(semi)stable representation, then there is a Zariski open subset of
$\theta$-(semi)stable representations in $rep_{\alpha}~Q$.

For dimension vectors $\alpha$ and $\beta$ we denote $\beta \rInto \alpha$ if there is an open
subset of $rep_{\alpha}~Q$ of representations containing a subrepresentation of dimension
vector $\beta$. Assume by induction we know for all $\gamma < \alpha$ the set of dimension
vectors $\delta$ such that $\delta \rInto \alpha$, then Schofield shows that
\[
\beta \rInto \alpha \quad \text{iff} \quad 0 = \underset{\gamma \rInto \beta}{max}~-\chi_Q(\gamma,\alpha-\beta)
\]
For given $\alpha$ let $S_{\alpha}$ be the finite set of dimension vectors $\beta$ such that
$\beta \rInto \alpha$, then
\begin{itemize}
\item{$\alpha$ is the dimension vector of a $\theta$-semistable representation of $Q$ if and
only if $\theta(\beta) \geq 0$ for all $\beta \in S_{\alpha}$.}
\item{$\alpha$ is the dimension vector of a $\theta$-stable representation of $Q$ if and only if
$\theta(\beta) > 0$ for all $\beta \in S_{\alpha}$.}
\end{itemize}

In many important applications it is easy to prove that certain low-dimensional dimension vectors
have $\theta$-stable representations. The following result can then be used to construct all
dimension vectors of $\theta$-stables lying in the positive integral span of these vectors.

\begin{theorem} \label{simples} Let $\beta_1,\hdots,\beta_l$ be dimension vectors of $\theta$-stable
representations of $Q$ and assume there are integers $m_1,\hdots,m_l \geq 0$ such that
\[
\alpha = m_1 \beta_1 + \hdots + m_l \beta_l \]
Then, $\alpha$ is the dimension vector of a $\theta$-stable representation of $Q$ if and only
if $\alpha' = (m_1,\hdots,m_l)$ is the dimension vector of a simple representation of the quiver $Q'$
on $l$ vertices $w_1,\hdots,w_l$ such that there are exactly
\[
\delta_{ij} - \chi_Q(\beta_i,\beta_j) \]
arrows from $w_i$ to $w_j$.
\end{theorem}

\begin{proof}
Let $W_i$ be a $\theta$-stable representation of $Q$ of dimension vector $\beta_i$ and consider
the $\alpha$-dimensional representation
\[
V = W_1^{\oplus m_1} \oplus \hdots \oplus W_l^{\oplus m_l} \]
It is clear from the definition that $(Q',\alpha')$ is the local quiver setting
corresponding to $V$.

As before, there is a semi-invariant $d_{\sigma}$ such that $V \in X_{\sigma}(\alpha) = rep_{\alpha}~\C Q_{\sigma}$
and $\C Q_{\sigma} \otimes V$ is a semi-simple representation of the universal
localization $\C Q_{\sigma}$.

If there are $\theta$-stable representations of dimension $\alpha$, then there is an open
subset of $X_{\sigma}(\alpha)$ consisting of $\theta$-stable representations. But we have seen
that they become simple representations of $\C Q_{\sigma}$. This means that every Zariski neighborhood
of $V \in rep_{\alpha}~\C Q_{\sigma}$ contains simple $\alpha$-dimensional representations. By
the \'etale local isomorphism there are $\alpha'$-dimensional simple representations of the quiver $Q'$.

Conversely, as any Zariski neighborhood of the zero representation in $rep_{\alpha'}~Q'$ contains
simple representations, then so does any neighborhood of $V \in rep_{\alpha}~\C Q_{\sigma} = X_{\sigma}(\alpha)$.
We have seen before that $X_{\sigma}(\alpha)$ consists of $\theta$-semistable representations and
that the $\theta$-stables correspond to the simple representations of $\C Q_{\sigma}$, whence
$Q$ has $\theta$-stable representations of dimension vector $\alpha$.
\end{proof}

We recall from \cite{LeBruynProcesi:1990} the classification of dimension vectors of simple
representations of quivers. $\alpha' = (m_1,\hdots,m_l)$ is the dimension vector of a simple
representation of $Q'$ if and only if the support of $\alpha'$ is a strongly connected quiver
(that is, every pair of vertices belongs to an oriented cycle) and we are in one
of the
following two cases
\begin{enumerate}
\item{Either $supp(\alpha') = \tilde{A_n}$ the extended Dynkin diagram with cyclic orientation for
some $n$ and $\alpha' \mid supp(\alpha') = (1,\hdots,1)$.}
\item{Or $supp(\alpha') \not= \tilde{A_n}$ and we have the following inequalities
\[
\chi_{Q'}(\epsilon_i,\alpha') \leq 0 \quad \text{and} \quad
\chi_{Q'}(\alpha',\epsilon_i) \leq 0 \]
for all vertices $w_i \in supp(\alpha')$ where $\epsilon_i = (\delta_{1i},\hdots,\delta_{li})$.}
\end{enumerate}
The above theorem asserts that the set of $\theta$-stable dimension vectors can be described
by a set of inequalities. For more results along similar lines we refer the reader to the
recent preprint \cite{DerksenWeyman:2000} of H. Derksen and J. Weyman.

In the following section we will apply these results to the classification of dimension vectors
corresponding to simple representations of the free product groups $\Z_p \ast \Z_q$.

\section{Torus knot groups.}

In this section we translate the problem of describing all irreducible representations of a
torus knot group to a quiver setting. Our method is based on earlier work of B. Westbury
\cite{Westbury}. In the next section we will apply the foregoing results to
determine the relevant dimension vectors.

Consider a solid cylinder $C$ with $q$ line segments on its curved face, equally spaced and parallel
to the axis. If the ends of $C$ are identified with a twist of
$2 \pi \tfrac{p}{q}$ where $p$ is an integer relatively prime to $q$, we obtain a single curve
$K_{p.q}$ on the surface of a solid torus $T$. If we assume that $T$ lies in $\R^3$ in the standard way,
the curve $K_{p.q}$ is called the $(p,q)$ torus knot.

The fundamental group of the complement $\R^3 - K_{p,q}$ is called the $(p.q)$-torus knot group
$G_{p,q}$ which has a presentation
\[
G_{p,q} = \pi_1(\R^3-K_{p,q}) \simeq \langle a,b~\mid~a^p=b^q \rangle \]
An important special case is $(p,q) = (2,3)$ in which case we obtain the three string braid group,
$G_{2,3} \simeq B_3$.

Recall that the center of $G_{p,q}$ is generated by $a^p$ and that the quotient group is the
free product of cyclic groups of order $p$ and $q$
\[
\overline{G_{p,q}} = \dfrac{G_{p,q}}{\langle~a^p~\rangle} \simeq \Z_p \ast \Z_q \]
As the center acts by scalar multiplication on any irreducible representation, the representation
theory of $G_{p,q}$ essentially reduces to that of $\Z_p \ast \Z_q$. Observe that in the
special case $(p,q) = (2,3)$ considered above, the quotient group is the modular group
$PSL_2(\Z) \simeq \Z_2 \ast \Z_3$.

Let $V$ be an $n$-dimensional representation of $\Z_p \ast \Z_q$, then the restriction of $V$ to
the cyclic subgroups $\Z_p$ and $\Z_q$ decomposes into eigenspaces
\[
\begin{cases}
V~\downarrow_{\Z_p} &\simeq S_1^{\oplus a_1} \oplus S_{\zeta}^{\oplus a_2} \oplus \hdots \oplus
S_{\zeta^{p-1}}^{\oplus a_p} \\
V~\downarrow_{\Z_q} &\simeq T_1^{\oplus b_1} \oplus T_{\xi}^{\oplus b_2} \oplus \hdots \oplus
T_{\xi^{q-1}}^{\oplus b_q}
\end{cases}
\]
where $\zeta$ (resp. $\xi$) is a primitive $q$-th (resp. $p$-th) root of unity and where $S_{\zeta^i}$
(resp. $T_{\xi^i}$) is the one-dimensional space $\C v$ with action $a.v = \zeta^i v$ (resp. $b.v =
\xi^i v$). Using these decompositions we define linear maps $\phi_{ij}$ as follows
\[
\begin{diagram}
S_{\zeta^{i-1}}^{\oplus a_i} & & \rTo^{\phi_{ij}} & & T_{\xi^{j-1}}^{\oplus b_j} \\
\dInto & & & & \uOnto \\
S_1^{\oplus a_1} \oplus S_{\zeta}^{\oplus a_2} \oplus \hdots \oplus
S_{\zeta^{p-1}}^{\oplus a_p} &\quad = & V & = \quad& T_1^{\oplus b_1} \oplus T_{\xi}^{\oplus b_2} \oplus \hdots \oplus
T_{\xi^{q-1}}^{\oplus b_q}
\end{diagram}
\]
This means that we can associate to an $n$-dimensional representation $V$ of $\Z_p \ast \Z_q$ a
representation of the full bipartite quiver on $p+q$ vertices
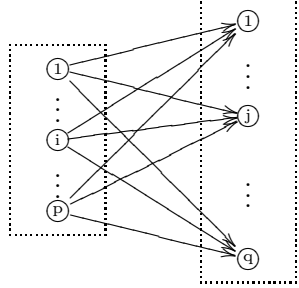
\[
\xy ;/r.15pc/:
(-10,-20); (-10,20) **@{.}; (10,20) **@{.}; (10,-20) **@{.}; (-10,-20) **@{.};
(30,-30); (30,30) **@{.}; (50,30) **@{.}; (50,-30) **@{.}; (30,-30) **@{.};

\POS (0,15) *\cir<4pt>{}*+{\txt\tiny{1}} ="v1"
   , (0,-15) *\cir<4pt>{}*+{\txt\tiny{p}} ="vk"
   , (0,8) *+{\vdots}
   , (0,-8) *+{\vdots}
   , (40,15) *+{\vdots}
   , (40,-10) *+{\vdots}
   , (0,0) *\cir<4pt>{}*+{\txt\tiny{i}} ="vj"
   , (40,25) *\cir<4pt>{}*+{\txt\tiny{1}} ="w1"
   , (40,5) *\cir<4pt>{}*+{\txt\tiny{j}} ="wi"
   , (40,-25) *\cir<4pt>{}*+{\txt\tiny{q}} ="wp"
\POS"v1" \ar "w1"
\POS"v1" \ar "wi"
\POS"v1" \ar "wp"
\POS"vk" \ar "w1"
\POS"vk" \ar "wi"
\POS"vk" \ar "wp"
\POS"vj" \ar "w1"
\POS"vj" \ar "wi"
\POS"vj" \ar "wp"
\endxy
\]
where we put at the left $i$-th vertex the space $S_{\zeta^{i-1}}^{\oplus a_i}$, on the right
$j$-th vertex the space $T_{\xi^{j-1}}^{\oplus b_j}$ and the morphism connecting the $i$-th
left vertex to the right $j$-vertex is the map $\phi_{ij}$. That is, to $V$ we associate a
representation $V_Q$ of dimension vector $\alpha = (a_1,\hdots,a_p;b_1,\hdots,b_q)$ and of course we
have that $a_1+\hdots+a_p = n = b_1 + \hdots + b_q$.

If $V$ and $W$ are isomorphic as $\Z_p \ast \Z_q$ representation, they have isomorphic weight
space decompositions for the restrictions to $\Z_p$ and $\Z_q$ and fixing bases in these weight
spaces gives isomorphic quiver representations $V_Q \simeq W_Q$. Further, observe that if $V$ is
a representation of $\Z_p \ast \Z_q$ then the matrix
\[
m(V_Q) = \begin{bmatrix} \phi_{11}(V_Q) & \hdots & \phi_{p1}(V_Q) \\
\vdots & & \vdots \\
\phi_{1q}(V_Q) & \hdots & \phi_{pq}(V_Q) \end{bmatrix} \]
is invertible and the determinant $d = det(m)$ is a semi-invariant for $rep_{\alpha}~Q$ of
weight $\theta = (\underbrace{-1,\hdots,-1}_p;\underbrace{1,\hdots,1}_q)$. That is, all representations of the form $V_Q$ are $\theta$-semistable representations
of $Q$.

\begin{lemma} If $V$ is an irreducible $\Z_p \ast \Z_q$ representation, then $V_Q$ is a $\theta$-stable
quiver representation. Conversely, if $W \in rep_{\alpha}~Q$ is $\theta$-semistable such that
$d(V) \not= 0$, then $W = V_Q$ for a representation $V$ of $\Z_p \ast \Z_q$ and if $W$ is $\theta$-stable
then $V$ is irreducible.
\end{lemma}

\begin{proof} Let $W \in rep_{\alpha}~Q$ such that $d(W) \not= 0$ and $U \rInto W$ a subrepresentation
of dimension vector $\beta = (c_1,\hdots,c_p;d_1,\hdots,d_q)$. Then, $\theta(U) \geq 0$ for if
$\theta(W) < 0$ then $\sum_i c_i > \sum_j d_j$ and the composed map $U_L \rTo^{\phi(U)} U_R$ must have a kernel
contradicting the fact that $\phi(V)$ is an isomorphism. If $\theta(U) = 0$, then $\sum_i c_i =
\sum_j d_j = m$ and the restriction of $\phi(V)$ to $U$ gives an linear isomorphism $U_L \simeq U_R$.
Finally, the decomposition $U_L$ determines the action of $\Z_p$ on an $m$-dimensional space $U'$ and
$U_R$ the action of $\Z_q$ on $U'$, that is, $U \simeq U'_Q$ for a representation $U'$ of
$\Z_p \ast \Z_q$.
\end{proof}

That is, the variety of semi-simple $n$-dimensional representations of $\Z_p \ast \Z_q$ decomposes
into components
\[
iss_n~\Z_p \ast \Z_q = \underset{\sum a_i=\sum b_j = n}{\bigsqcup}~iss_{\alpha}~\C Q_{\sigma} \rInto
\underset{\sum a_i = \sum b_j = n}{\bigsqcup}~M^{ss}_{\alpha}(Q,\theta) \]
(where $\sigma$ is the map $\phi$ described before)
which are affine open subsets of the projective moduli spaces of $\theta$-semistable representations
of $Q$. In the next section we will determine the dimension vectors $\alpha$ corresponding to
irreducible representations.

As an example, consider the infinite dihedral group $D_{\infty} = \Z_2 \ast \Z_2$, then the
corresponding quiver is $Q=\tilde{A_4}$ with orientation
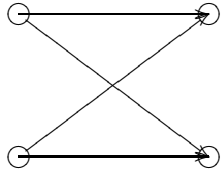
\[
\xy ;/r.15pc/:
\POS (0,15) *\cir<4pt>{}*+{\txt\tiny{}} ="v1"
   , (0,-15) *\cir<4pt>{}*+{\txt\tiny{}} ="v2"
   , (40,15) *\cir<4pt>{}*+{\txt\tiny{}} ="w1"
   , (40,-15) *\cir<4pt>{}*+{\txt\tiny{}} ="w2"
\POS"v1" \ar "w1"
\POS"v1" \ar "w2"
\POS"v2" \ar "w1"
\POS"v2" \ar "w2"
\endxy
\]
For $\alpha = (1,1;1,1)$ we have $M^{ss}_{\alpha}(Q,\theta) \simeq \PP^1$ and
$iss_{\alpha}~D_{\infty} = \Af^1$, the one point missing is the representation $V$ where
all maps are the identity because in that case $d(V) = 0$.

\section{Dimension vectors of irreducible representations.}

In this section we will characterize the dimension vectors $\alpha$ such that $\Z_p \ast \Z_q$
has irreducible representations $V$ such that $V_Q$ has dimension vector $\alpha$. In view of the
above we have to determine the dimension vectors of $\theta$-stable representations of $Q$.
There are some obvious $1$-dimensional irreducible representations of $\Z_p \ast \Z_q$
\[
V_{ij} = \C v \quad \text{with $a.v = \zeta^{i-1} v$ and $b.v = \xi^{j-1} v$.} \]
which have dimension vector $\alpha_{ij} = (\delta_{1i},\hdots,\delta_{pi};\delta_{1j},\hdots,
\delta_{qj})$.

We consider the quiver $\Gamma$ on $i.j$ vertices $v_{ij}$ (corresponding to the $\theta$-stable representations
$V_{ij}$) such that the number of arrows from $v_{ij}$ to $v_{kl}$ is equal to
\[
\delta_{ij,kl} - \chi_Q(\alpha_{ij},\alpha_{kl}) \]
Given the special form of the full bipartite quiver $Q$ it is easy to verify that
\[
\#~\{~a \in \Gamma_a~\mid~\xy
\POS (15,0) *\cir<7pt>{}*+{\txt\tiny{$v_{ij}$}} ="v2"
   , (0,0) *\cir<7pt>{}*+{\txt\tiny{$v_{kl}$}} ="v3"
\POS"v2" \ar "v3"_{a} \endxy~\} = \begin{cases}
1 \quad &\text{if $i \not= k$ and $j \not= l$} \\
0 \quad &\text{otherwise.}
\end{cases}
\]
For example, in the modular case $PSL_2(\Z) = \Z_2 \ast \Z_3$ the quiver $\Gamma$ has the form
of a large cycle with subsequent vertices $v_{21},v_{12},v_{23},v_{11},
v_{22},v_{13}$ and only one arrow in each direction between two consecutive
vertices.
%\[
%\xy
%\POS (0,0) *\cir<6pt>{}*+{\txt\tiny{$v_{11}$}} ="S11",
%\POS (14,14) *\cir<6pt>{}*+{\txt\tiny{$v_{23}$}} ="S23",
%\POS (-14,14) *\cir<6pt>{}*+{\txt\tiny{$v_{22}$}} ="S22",
%\POS (14,34) *\cir<6pt>{}*+{\txt\tiny{$v_{12}$}} ="S12",
%\POS (-14,34) *\cir<6pt>{}*+{\txt\tiny{$v_{13}$}} ="S13",
%\POS (0,48) *\cir<6pt>{}*+{\txt\tiny{$v_{21}$}} ="S21",
%\POS"S11" \ar@/^2ex/ "S23"
%\POS"S23" \ar@/^2ex/ "S11"
%\POS"S12" \ar@/^2ex/ "S23"
%\POS"S23" \ar@/^2ex/ "S12"
%\POS"S12" \ar@/^2ex/ "S21"
%\POS"S21" \ar@/^2ex/ "S12"
%\POS"S21" \ar@/^2ex/ "S13"
%\POS"S13" \ar@/^2ex/ "S21"
%\POS"S13" \ar@/^2ex/ "S22"
%\POS"S22" \ar@/^2ex/ "S13"
%\POS"S11" \ar@/^2ex/ "S22"
%\POS"S22" \ar@/^2ex/ "S11"
%\endxy
%\]
We consider the direct sum of $\theta$-stable representations
\[
V = \oplus_{i,j}~V_{ij}^{\oplus m_{ij}} \]
that is, the dimension vector of $V$ is $\alpha = \sum_{i,j} m_{ij} \alpha_{ij}$. In view of
theorem~\ref{simples} to verify that there are $\theta$-stable representations of dimension vector
$\alpha$ is equivalent to verifying that $\gamma = (m_{11},\hdots,m_{pq})$ is the dimension vector
of a simple representation of $\Gamma$.

\begin{theorem} $\alpha = (a_1,\hdots,a_p;b_1,\hdots,b_q)$ with $\sum_i a_i = n = \sum_j b_j$ 
is the dimension vector of a 
$\theta$-stable representation of $Q$ if and only if
\[
n = 1 \quad \text{or} \quad a_i + b_j \leq n \]
for all $1 \leq i \leq p$ and $1 \leq j \leq q$.
\end{theorem}

\begin{proof}
Assume that $\alpha$ is the dimension vector of a $\theta$-stable representation of $Q$, then
the set of irreducible representations of $\Z_p \ast \Z_q$ is an open subset of $M^{ss}_{\alpha}(Q,\theta)$.
Assume that $a_i + b_j > n$ and consider the vertex spaces (the eigenspaces) $S_i^{\oplus a_i}$ and $T_j^{\oplus b_j}$ as
subspaces of the $n$-dimensional representation $V' = \C^n$ (via the identification given by the matrix
$m$). Then $W = S_i^{\oplus a_i} \cap T_j^{\oplus b_j} \not= 0$ is a subrepresentation of $V$ of
dimension vector $a \alpha_{ij}$ for some $a \in \N_+$. But $\theta(W) = 0$ so $V'$ cannot be $\theta$-stable,
a contradiction.

Conversely, assume the numerical condition is satisfied and consider the direct sum of $\theta$-stable
representations
\[
V = \oplus_{i,j} V_{ij}^{\oplus m_{ij}} \]
We note that
     $$a_{i}=\sum_{j=1}^q m_{ij}\qquad \text{ and }\qquad
     b_{l}=\sum_{i=1}^p m_{il}.$$
Let $I_{p}\in M_{p}(\C)$ be the $p\times p$ identity matrix and let
     $A_{p}\in M_{p}(\C)$ be
     the $p\times p$ matrix of the form
     $$A_{p}=
     \begin{pmatrix}
	0 & -1 & \hdots & -1 & -1\\
	-1 & 0 & \hdots & -1 & -1\\
	\vdots & \vdots & \ddots & \vdots & \vdots\\
	-1 & -1 & \hdots & 0 & -1\\
	-1 & -1 & \hdots & -1 & 0
     \end{pmatrix}.
     $$
     Then the Euler form of the local quiver $\Gamma$ is the symmetric matrix 
     $$\chi_{\Gamma}=
     \begin{pmatrix}
	I_p & A_p & \hdots & A_p & A_p\\
	A_p & I_p & \hdots & A_p & A_p\\
	\vdots & \vdots & \ddots & \vdots & \vdots\\
	A_p & A_p & \hdots & I_p & A_p\\
	A_p & A_p & \hdots & A_p & I_p
     \end{pmatrix}\qquad \in \qquad M_{q}(M_{p}(\C)).
     $$
When $n=1$, we have that $V=V_{11}$ is obviously a $\theta$-stable
     representation.
When $n=2$, we notice that $\gamma$ is the dimension
     vector of a simple representation if and only if
     $$\gamma \arrowvert supp(\gamma)=(1,1;1,1)\qquad \text{ because }\qquad
     supp(\gamma)=\widetilde{A}_{2}$$
Now, consider $n \geq 3$ and consider the dimension vector
$\gamma = (m_{11},\hdots,m_{pq})$. We have to verify that $\gamma$ is the dimension vector
of a simple representation of $\Gamma$ which, by symmetry of $\chi_{\Gamma}$, amounts to
checking that
     $$
     \chi_{\Gamma}(\gamma,\epsilon_{kl})=
     \chi_{\Gamma}(\epsilon_{kl},\gamma)\leq 0$$
     for all $1\leq k\leq p$ and $1\leq l\leq q$ where $\epsilon_{kl} = (\delta_{ij,kl})$ are the
     standard base vectors. Computing the left hand term this is equivalent to
     $$m_{kl}+\sum_{k\neq i =1}^p \sum_{l\neq j=1}^q -m_{ij}\leq 0.$$
     or
     $$m_{kl}+\sum_{k\neq i=1}^p m_{il}\leq
     \sum_{k\neq i =1}^p \sum_{l\neq j=1}^q m_{ij}+\sum_{k\neq i=1}^p m_{il}$$
Resubstituting the values of $a_i$ and $b_l$ in this expression we see that this is
equivalent to
     $$b_{l}\leq \sum_{k\neq i=1}^p a_{i}=n-a_{k}$$
Therefore, the condition is satisfied if for all 1$\leq k\leq p$ and $1\leq l\leq q$ we have
     $$a_{k}+b_{l}\leq n.$$
\end{proof}

In the special case of $PSL_2(\Z) = \Z_2 \ast \Z_3$, our condition on the
dimension vector $\alpha = (a_1,a_2;b_1,b_2,b_3)$ is equivalent to
\[
a_i + b_j \leq n = a_1+a_2 \quad \text{whence} \quad b_j \leq a_i \]
for all $1 \leq i \leq 2$ and $1 \leq j \leq 3$ which was the criterium found by B. Westbury in
\cite{Westbury}.

Using the results of the present paper, a complete characterization is given
in \cite{ABV:2000} of the dimension vectors $\alpha$ for
which the moduli space $M^{ss}_{\alpha}(Q,\theta)$ is smooth.

\end{document}